\documentclass[twoside,11pt,leqno]{article}
\newcommand{\oR}{{\mathbb R}}

\newcommand{\oN}{{\mathbb N}}
\newcommand{\xx}{\mathbf{x}}
\newcommand{\yy}{\mathbf{y}}
\newcommand{\zz}{\mathbf{z}}

\usepackage{amssymb}
\newtheorem{defn}{Definition}
\newtheorem{res}{Lemma}

\newtheorem{thm}{Theorem}

\newtheorem{ex}{Example}
\newenvironment{proof}[1]{\noindent{\bf Proof of #1:} }{\hfill $\square$ \\}

\begin{document}
\bibliographystyle{plain}
\thispagestyle{empty}
\begin{center}
{\Large\bf Nearest-neighbour Markov point processes on \\
graphs with Euclidean edges}\\[.4in] 

\noindent
{\large M.N.M. van Lieshout}\\[.1in]
\noindent
{\em CWI, P.O.~Box 94079, NL-1090 GB  Amsterdam\\
University of Twente, P.O.~Box 217, NL-7500 AE Enschede\\
The Netherlands}
\end{center}

\begin{verse}
{\footnotesize
\noindent
{\bf Abstract}\\
\noindent
We define nearest-neighbour point processes on graphs with Euclidean
edges and linear networks. They can be seen as the analogues of renewal 
processes on the real line. We show that the Delaunay neighbourhood 
relation on a tree satisfies the Baddeley--M\o ller consistency conditions
and provide a characterisation of Markov functions with respect to this
relation. We show that a modified relation defined in terms of the local 
geometry of the graph satisfies the consistency conditions for all graphs 
with Euclidean edges. \\[0.1in]
\noindent
{\em Keywords \& Phrases:} Delaunay neighbours; Graph with Euclidean edges;
Linear network; Markov point process; Nearest-neighbour interaction; 
Renewal process. \\[0.1in]
\noindent
{\em 2010 Mathematics Subject Classification:}
60G55, % point processes     
60D05. % stochastic geometry 
}
\end{verse}

\section{Introduction}

In recent years, a theory of point processes on linear networks has
been emerging so as to be able to analyse, for example, the prevalence 
of accidents on motorways, the occurrence of street crimes and other
data described in the first chapter of the pioneering monograph by 
Okabe and Sugihara \cite{OkabSugi12}. Although there exists a mature
theoretical framework for point processes on Euclidean spaces 
\cite{DaleVere88}, the development of a similar theory on linear  
networks is complicated by the geometry inherent in the network. In 
particular, it is not possible to define strictly stationary models, 
as the network may not be closed under translations. For this reason, 
most attention has focussed on the development of second order summary 
statistics \cite{Raksetal17}.

Little attention has been paid to model building with a few notable 
exceptions. The first serious work in this direction seems to be that
by Baddeley et al.~\cite{Baddetal17}, who construct certain types
of Cox processes as well as a Switzer-type and a cell process. The
authors conclude that familiar procedures for constructing models 
tend not to produce processes on a linear network that are pseudostationary
with respect to the shortest path distance, except when the network 
is a tree -- an unrealistic assumption for a road network. Another 
important contribution is the work by Anderes et al.~\cite{Moll16} 
who expand the modelling framework in various directions. They relax 
the assumption of \cite{Baddetal17,OkabSugi12,Raksetal17}
that a linear network consists of a finite union of straight line 
segments that intersect only in vertices, in the sense that the 
segments are replaced by parametrised rectifiable curves that
may or may not overlap. The parametrisations have the additional 
advantage of naturally defining a weighted shortest path distance. 
In the motivating example where the linear network represents a 
road network, such a generalisation allows for bridges or tunnels 
and for distance to be measured in travel time where appropriate. 
Additionally, Anderes et al.~\cite{Moll16} construct log-Gaussian 
Cox processes in terms of a Gaussian process on the network specified 
by an isotropic covariance function. 

All models discussed so far are clustered in nature, that is, exhibit
a positive association between the points. In this paper, our aim is
to develop appropriate analogues of renewal processes, exploiting 
the one-dimensional nature of a linear network. Recall that renewal
processes exhibit the property that the probability (in an infinitesimal
sense) of an event at a given location conditional on the realisation
of the process elsewhere depends only on the two nearest points, 
regardless of how far away they may be. Such models are known as 
nearest-neighbour Markov point processes \cite{BaddMoll89} and may 
incorporate both inhibition and clustering \cite{Lies00}. In contrast to 
Cox models, the second order summary statistics may not be available
in closed form, but the conditional intensities and likelihood are.
The latter can be expressed as a product of interaction functions, 
which may be chosen to be isotropic.

The plan of this paper is as follows. In Section~\ref{S:prelim}, 
we recall the definitions of Anderes et al.~\cite{Moll16} regarding 
graphs with Euclidean edges, the weighted shortest path metric 
thereon and Poisson process defined on them. In Section~\ref{S:nnMp}, 
we extend the notion of a Markov point process with respect to the
Delaunay nearest-neighbour relation \cite{BaddMoll89} to graphs
with Euclidean edges and state our main results. More specifically,
we show that the Delaunay relation on a tree satisfies the 
Baddeley--M\o ller consistency conditions and we provide a 
characterisation of the Markov functions with respect to this 
relation. We then use the graph structure to define a modified 
Delaunay relation and show that it satisfies the consistency 
conditions on any graph with Euclidean edges. The proofs are 
given in Sections~\ref{S:proofs}--\ref{S:theorems}.

\section{Preliminaries} \label{S:prelim}

\subsection{Graphs with Euclidean edges}

In their pioneering monograph on the subject, Okabe and Sugihara 
\cite[page 31]{OkabSugi12} define a network as a finite union
\[
L = \bigcup_{i=1}^n l_i, \quad n \in \oN,
\]
of straight line segments $l_i$ in $\oR^2$ or $\oR^3$ that intersect
only at their endpoints in such a way that $L$ is connected. The 
representation is not unique since a line segment may arbitrarily
be split in two pieces without affecting the union $L$. 

A more general definition is given by Anderes, M\o ller and Rasmussen
\cite{Moll16}. They replace the straight line segments by curves that
are parametrised by bijections. In order to define a weighted shortest 
path distance on the graph, we impose the further condition that these 
parametrisations are diffeomorphisms and follow \cite{Moll16} to define 
the weight of an edge as the `length' induced by the parametrisation.

\begin{defn} \label{d:graph}
A graph with Euclidean edges in $\oR^2$ is a triple 
$(V, E = (e_i)_i, \Phi = (\phi_i)_i)$ such that
\begin{itemize}
\item $(V,E)$ is a finite, simple connected graph, i.e. has neither
loops nor multiple edges;
\item every edge $e_i = \{ v_i^1, v_i^2 \} \in E$, 
$v_i^1, v_i^2 \in V$, is parametrised by the inverse of a 
bijection $\phi_i$, that is, $\phi_i^{-1}: J_i \to \oR^2$ for a non-empty 
open interval $J_i \subset \oR$ with endpoints $\phi_i(v_i^j)$, $j=1, 2$.
\end{itemize}
\end{defn}

In other words, an edge is associated with a set 
$\phi_i^{-1}(J_i) \subset \oR^2$ that does not contain the endpoints. 

A graph with Euclidean edges gives rise to the space
\[
L = (\{ 0 \} \times V ) \cup
  \bigcup_{i=1}^{n(E)} ( \{ i \} \times \phi_i^{-1}(J_i),
\]
where $n(E) < \infty$ is the cardinality of $E$. The labels $i$ serve 
to identify the edges 
and will prevent paths from `jumping from one edge to another' 
in case their interiors overlap. For instance, if $L$ represents a 
road network, overlap is typically  present due to tunnels and bridges
\cite{Moll16}.

As an aside, if there is no overlap between the edge interiors, one may 
drop the labels and simply consider the disjoint union
\[
V \cup \bigcup_{i=1}^{n(E)} \phi_i^{-1}(J_i),
\]
which in turn reduces to the classic linear networks if all edges
$\phi_i^{-1}(J_i)$ are straight line segments. From now on, we 
will work with the general space $L$ including the labels.

\subsection{Weighted shortest path metric}

The family $\Phi$ of parametrisations that is part of the definition
of a graph with Euclidean edges can be used to define concepts of length
and distance \cite{Moll16}.

\begin{defn} \label{d:metric}
Let $(V, E, \Phi)$ be a graph with Euclidean edges. For every 
$i=1, \dots, n(E)$, define the $\phi_i^{-1}$-induced length measure 
on the $\sigma$-algebra
\(
\{ \phi_i^{-1}(B): B \subset \bar J_i \mbox{ Borel} \}
\) 
on $\phi_i^{-1}(\bar J_i)$ generated by the functions $\phi_i$ as follows.
The induced length of the set $\phi_i^{-1}(B)$ is equal to the Euclidean 
length of $B$ in the closed interval $\bar J_i$. In particular, the edge 
$e_i = \{ v_i^1, v_i^2 \}$ has length $|\phi(v_i^1) - \phi(v_i^2)|$ equal 
to the Euclidean length of $J_i$.
\end{defn}

If we use the arc length parametrisation, the induced length of 
Definition~\ref{d:metric} corresponds to the usual Euclidean length. A 
different set of parametrisations may induce a different length measure. 

\begin{defn} \label{d:walk}
Let $(V, E, \Phi)$ be a graph with Euclidean edges. 
A walk between two elements $(i,x)$ and $(j,y)$ of $L$ travels 
alternatingly from $(i,x)$ to $(j,y)$ along a finite number of nodes 
and edges such that the nodes are endpoints of the edges. If $i\neq 0$, 
the first bit of the walk travels from $x$ along $\phi_i^{-1}( J_i)$; 
similarly for $j\neq 0$, the last bit of the walk is along 
$\phi_j^{-1}( J_j)$  to $y$. In particular, a walk between two points 
along the same curve $\phi_i^{-1}(J_i)$ does not reverse its tracks.
\end{defn}

\begin{defn} \label{d:path}
Let $(V, E, \Phi)$ be a graph with Euclidean edges. 
A path between two different elements $(i,x)$ and $(j,y)$ of $L$
is a walk in which all edge segments and all vertices are different. 
The weight of the path is the sum of the lengths of the edge segments 
in it, and the shortest path distance $d_G( (i,x), (j,y) )$ between 
two elements $(i,x)$ and $(j, y)$ of $L$ is the smallest weight carried 
by any path between them. 
\end{defn}

Note that it is possible that the shortest weighted path between two 
vertices joined by an edge is not along that edge! The assumption that 
the parametrisations are diffeomorphisms ensures that the edge lengths 
are finite. Moreover, the weighted shortest path distance on a graph with 
Euclidean edges defines a metric. For further details, see \cite{Moll16}.

\subsection{Poisson process on graphs with Euclidean edges}

Our goal in this paper is to define point process on graphs with
Euclidean edges in terms of a density with respect to a unit rate
Poisson process \cite{Lies00}. To do so, we recall Anderes et al.'s
definition of a Poisson process on a graph with Euclidean edges 
\cite{Moll16}.

\begin{defn} \label{d:Lebesgue}
Let $(V, E = (e_i)_i, \Phi = (\phi_i)_i)$ be a graph with Euclidean
edges and $L$ the corresponding network. The Lebesgue--Stieltjes
measure $\lambda_G$ is defined as follows. For every $i$ and every 
set $B_i$ in the $\sigma$-algebra generated by $\phi_i$, set
\[
\lambda_G^i(B_i) = \int_{J_i} 1\{ \phi_i^{-1}(t) \in B_i \} 
\left| \frac{d}{dt} \phi_i^{-1}(t) \right| \, dt,
\]
where $| \cdot |$ denotes the norm of the gradient. Then
\[
\lambda_G\left( \cup_{i=1}^{n(E)} ( \{ i \} \times B_i ) \right) = 
\sum_{i=1}^{n(E)} \lambda_G^i( B_i )
\]
is a measure on $L$ equipped with the finite disjoint union $\sigma$-algebra
in which a set is measurable if and only if it can be written as
$\cup_{i=0}^{n(E)} ( \{ i \} \times B_i )$ with $B_i$ in the 
$\sigma$-algebra generated by $\phi_i$ for $i>0$ and the power set 
of $V$ for $i=0$. By default, the set $\{ 0 \} \times V$ has 
Lebesgue--Stieltjes measure zero.
\end{defn}

Note that Definition~\ref{d:Lebesgue} does not depend on the 
parametrisations. 

\begin{defn} \label{d:Poisson}
Let $(V, E = (e_i)_i, \Phi = (\phi_i)_i)$ be a graph with Euclidean
edges and $L$ the corresponding network. The unit rate Poisson process 
on $L$ is defined as follows. For every $i$ and every set $B_i$ in the 
$\sigma$-algebra generated by $\phi_i$,
\begin{itemize}
\item the number of points in $\{ i \} \times B_i$ is Poisson distributed 
with expectation $\lambda_G( \{ i \} \times B_i )$;
\item given that the number of points falling in $\{ i \} \times B_i$ is 
$n_i$, these $n_i$ points are independent and identically distributed according 
to the probability density function $1/\lambda_G( \{ i \} \times B_i )$.
\end{itemize}
\end{defn}

In words, the unit rate Poisson process scatters a Poisson number of 
points independently and uniformly on every edge, and the average number 
of such points is equal to the arc length of the edge. The definition does 
not depend on the parametrisation.

The integral of a measurable function $f: \oR^2 \to \oR^+$ with respect to 
the Lebesgue--Stieltjes measure $\lambda_G$ is defined as the sum of the
line integrals of $f$ along the rectifiable curves parametrised by the 
diffeomorphism $\phi_i^{-1}: J_i \to \oR^2$:
\[
\int_L f d\lambda_G = \sum_{i=1}^{n(E)}
\int_{J_i} f(\phi_i^{-1}(t)) \left| \frac{d}{dt} \phi_i^{-1} (t) \right| \, dt.
\]
Higher order integrals are defined analogously. Note that the definition
does not depend on the parametrisation.

\begin{ex}
As an example, consider the functions $f_G^i$ defined on $\phi_i^{-1}(J_i)$ by 
\[
f_G^i(x) = \left | { (\frac{d}{dt}\phi_i^{-1}) (\phi_i(x) )} \right|^{-1} .
\]
Then, for $B_i$ in the $\sigma$-algebra generated by $\phi_i$, 
\begin{eqnarray*}
\int_L 1 \{ x \in B_i \} \, f_G^i(x) \, d\lambda_G(x) & = &
\int_{J_i} 1\{ \phi_i^{-1}(t) \in B_i \} \left|
{ (\frac{d}{dt}\phi_i^{-1}) ( \phi_i( \phi_i^{-1}(t) ) } \right|^{-1} \, 
\left| \frac{d}{dt} \phi_i^{-1}(t) \right| \, dt \\
& =  & | \phi_i(B_i) |,
\end{eqnarray*}
the Euclidean length of $\phi_i(B_i)$. In other words, the functions $f_G^i$,
$i=1, \dots, n(E)$, define the weighted shortest path distance $d_G$ on $L$. 
\end{ex}

A point process $X$ on $L$ is said to have probability density $p$ with 
respect to the unit rate Poisson process if 
\begin{equation}
\label{e:dens}
P( X \in F) = \sum_{n=0}^\infty \frac{e^{-\lambda_G(L)}}{n!}
\int_L \cdots \int_L 1\{  \{ x_1, \dots, x_n \} \in F \} \,
p( \{ x_1, \dots, x_n \} ) \,
 d\lambda_G(x_1) \cdots d\lambda_G(x_n)
\end{equation}
for all $F$ in the usual $\sigma$-algebra on locally finite point
configurations in $L$ generated by the counts $X(B)$ with $B$ in the
finite disjoint union $\sigma$-algebra on $L$ \cite{DaleVere88}. For 
convenience's sake, we will use the notation $\xx = \{ x_1, \dots, x_n \}$ 
for configurations of finitely many distinct points in $L$.
For further details on graphs with Euclidean edges, the reader may
consult Anderes et al.\ \cite{Moll16}.

\section{Nearest-neighbour Markov point processes}
\label{S:nnMp}

The purpose of this section is to define appropriate analogues of the 
well-known class of renewal processes on the real line. Intuitively speaking,
we are looking for a class of point processes in which the conditional
behaviour at a given point depends on the remainder of the configuration
only through its `nearest neighbours'. We shall show that in the case 
that the network $L$ is a tree, the shortest path distance may be used 
to define which points are each other's nearest neighbours. In the general 
case, we will use the geometry of the network to define a local 
neighbourhood relation. 

\subsection{The Delaunay relation}

In this section, we adapt Baddeley and M\o ller's definition of configuration 
dependent neighbourhoods and cliques in Euclidean spaces \cite{BaddMoll89} 
to our context.

\begin{defn} \label{d:Delaunay}
Let $(V, E , \Phi)$ be a graph with Euclidean edges and $L$ the corresponding 
network. Let $\xx \subseteq L$ be a finite configuration of distinct points 
and define the Delaunay relation $\sim_{\xx}$ as the symmetric, reflexive 
relation on $\xx$ given by 
\[
x_1 \sim_{\xx} x_2 \Leftrightarrow  
C( x_1 | \xx ) \cap C( x_2 | \xx ) \neq \emptyset
\]
where, for $x_j \in \xx$, $j=1, 2$,
\[
C( x_j | \xx ) = \{ (i,y) \in L : d_G( (i,y), x_j ) \leq d_G( (i,y), x )
\mbox{ for all } x \in \xx \}
\]
is the Voronoi cell of $x_j$ in $\xx$. 
The $\xx$--neighbourhood of a subset $\zz \subseteq \xx$ is 
defined as 
\[
N(\zz | \xx ) = \{ x\in \xx : 
x \sim_{\xx} z \mbox{ for some } z \in \zz \}.
\]
The configuration $\zz$ is an $\xx$-clique if for each $z_1, z_2 \in\zz$,
$z_1 \sim_{\xx} z_2$. By convention, the empty set and singletons are
cliques too.
\end{defn}

A few remarks are in order. Although the space $L$ can be seen as
one dimensional, the Voronoi cells are not necessarily line segments
and the number of neighbours a point can have is not restricted to two.
Some further elementary properties of the relation are collected in the 
next lemma.

\begin{res}
\label{L:hereditary}
Let $(V, E , \Phi)$ be a graph with Euclidean edges and $L$ the corresponding 
network. The Delaunay relation on $L$ satisfies the following properties.
\begin{itemize}
\item If $\chi(\yy|\xx) = 1$ then also $\chi(\yy|\zz) = 1$ for all
$\yy \subseteq \zz \subseteq \xx$.
\item If $\chi(\yy|\zz) = 0$ then also $\chi(\yy|\xx) = 0$ for all
$\yy \subseteq \zz \subseteq \xx$.
\end{itemize}
Here $\chi(\yy|\xx) = 1$ if $\yy$ is an $\xx$-clique and zero otherwise.
\end{res}

The proof can be found in Section~\ref{S:proofs}.

\subsection{The Delaunay relation on trees}

Recall that $(V,E)$ is said to be a tree if it has no cycles, that is,
there is no closed path $(v_0, v_1, \dots, v_{p-1}, v_0)$, $v_i \in V$,
of positive length ($p > 0$). It is well-known that a graph is a tree 
if and only if there is exactly one path between any two vertices 
\cite{BondMurt08}.

A graph with Euclidean vertices $(V, E, \Phi)$ is said to be a tree if 
$(V,E)$ is. 

\begin{res}
\label{L:unique}
A graph with Euclidean edges is a tree if and only if there is exactly
one path between any two points $(i,x)$ and $(j,y)$ in $L$.
\end{res}

We are particularly interested in the geometrical arrangement of the paths 
between three points.

\begin{res}
\label{L:wheel}
Let $(V, E, \Phi)$ be a graph with Euclidean edges that is a tree
and $L$ its associated network. Consider a triple 
$\yy = \{ y_1, y_2, y_3 \} \subseteq L$. Then there exist unique paths 
between the elements of $\yy$ that 
\begin{itemize}
\item either form a wheel with three stokes of strictly positive length 
emanating from a vertex $(0,v)\in L$,
\item or combine into a single path.
\end{itemize}
\end{res}

The proofs of Lemma~\ref{L:unique} and \ref{L:wheel} will be given in 
Section~\ref{S:proofs}. 

\

A configuration $\yy\subseteq L$ is said to be in general position
if no three points lie on the boundary of the same $d_G$-ball. In 
other words, by Lemma~\ref{L:wheel}, a configuration in general position 
cannot form a wheel with three or more equally long spokes. Such a 
configuration would lead to Delaunay cliques with more than two elements.
Restriction to configurations in general position ensures that the 
size of cliques is at most two. Clearly, the class of all configurations 
in general position is hereditary.

\begin{res}
\label{L:clique}
Let $(V, E, \Phi)$ be a graph with Euclidean edges that is a tree
and $L$ its associated network. Then the clique sizes with respect 
to the Delaunay relation are at most two on the class of configurations 
in general position. Moreover, for all $\yy = \{ y_1, y_2 \} \subseteq \xx$ 
with $\xx$ in general position,
\(
\chi(\yy | \xx) = 1 
\)
if and only if the midpoint of $\yy$ with respect to the weighted
shortest path metric along the unique path between $y_1$ and
$y_2$ lies in $C(y_1 | \xx ) \cap C(y_2 | \xx)$.
\end{res}

Since the clique structure depends on the configuration, 
consistency conditions must be imposed on the family of 
neighbourhood relations \cite{BaddMoll89}. 

\begin{defn} \label{d:consist}
Let $(V, E, \Phi)$ be a graph with Euclidean edges and $L$ its 
associated network. Consider finite point configurations 
$\yy \subseteq \zz \subset L$ and points $u,v\in L$ such that 
$u,v\not\in\zz$. Then the Baddeley--M\o ller consistency conditions 
read as follows:
\begin{itemize}
\item[(C1)] $\chi(\yy | \zz)\neq \chi(\yy | \zz\cup\{ u\})$ 
implies $\yy\subseteq N(\{ u\} | \zz\cup\{ u\}) $;
\item[(C2)] if $u\not \sim_{\xx} v$ for $\xx = \zz \cup\{ u, v \}$, 
then 
\[
    \chi(\yy | \zz\cup\{ u\}) +  \chi(\yy | \zz\cup\{ v\}) = 
    \chi(\yy | \zz ) + \chi(\yy | \xx)
\]
where $\chi$ is the clique indicator function, i.e.\ 
$\chi(\yy | \xx) = 1\{ \yy$ is an $\xx$-clique$\}$.
\end{itemize}
\end{defn}

In general, the consistency conditions of Definition~\ref{d:consist}
do not hold, as illustrated by the following counterexample.

\begin{ex} \label{e:triangle}
Take $V = \{ (-1, 0), (1, 0), (0, 1) \} \subseteq \oR^2$ and connect all 
three vertices by straight line segments to form the complete graph 
$(V, E)$ and the linear network $L$. Let $\zz$ be a configuration
containing a single point on each edge. Then  $\chi( \zz | \zz ) = 1$. 
However, any additional point $u\not \in \zz$ would split up two of 
the points in $\zz$. Hence $\chi(\zz | \zz \cup \{ u \}) = 0$, 
although the third point of $\zz$ is no Delaunay neighbour of $u$ in
$\zz\cup\{ u \}$. Consequently, (C1) does not hold. Upon the addition 
of a second point $v$ that is not a Delaunay neighbour of $u$ in the 
resulting five point configuration, 
\(
\chi(\zz | \zz \cup \{ u, v\}) = 
\chi(\zz | \zz \cup \{ u \}) = \chi(\zz | \zz \cup \{ v \})
 = 0
\)
although $\chi(\zz | \zz) = 1$, in violation with (C2).
\end{ex}

The main theorem of this section is the following. Its proof will be
given in Section~\ref{S:theorems}.

\begin{thm} \label{t:Delaunay}
Let $(V, E, \Phi)$ be a graph with Euclidean edges that is a tree
and $L$ its associated network. Then the Delaunay relation satisfies 
(C1)--(C2) on the family of configurations in general position.
\end{thm}

\subsection{Markov functions}

We are now ready to define Markov functions on graphs with 
Euclidean edges, in analogy with the spatial models of \cite{BaddMoll89}.

\begin{defn} \label{d:nnMpp}
Let $(V, E, \Phi)$ be a graph with Euclidean edges and $L$ the 
corresponding network. Let $\sim_{\xx}$  be a family of reflexive, 
symmetric relations on finite point configurations $\xx$ in $L$. 
Then a function $p$ from the set of finite point configurations into 
$[0, \infty)$ is a  {\em Markov function\/} with respect to $\sim_\xx$ 
if for all $\xx$ such that $p(\xx)>0$,
\begin{description}
\item[(a)]  $p(\yy) > 0$ for all
$\yy \subseteq \xx$;
\item[(b)] for all $u\in L$, ${p(\xx \cup \{u\})} / {p(\xx)}$
depends only on $u$,
on $N( \{ u \} | \xx \cup\{ u\}) \cap \xx =
\{x \in \xx : x \sim_{\xx \cup\{ u\} } u \}$
and on the relations $\sim_{\xx} , \sim_{\xx \cup \{ u \} }$
restricted to $N(\{ u \} | \xx \cup \{ u \} )$.
\end{description}
\end{defn}

The next theorem provides a Hammersley--Clifford factorisation.
Similar results for spatial point processes in Euclidean spaces
can be found in \cite{BaddMoll89}, \cite{Baddetal96}, 
\cite{HaggLiesMoll99} or \cite{RiplKell77}. 
Recall that a function $\gamma$ from the space of finite point 
configurations into $[0,\infty)$ is an interaction function
\cite{BaddMoll89} if the following properties hold. If 
$\gamma(\xx) > 0$ then also 
(a) $\gamma(\yy) > 0$ for all $\yy\subseteq\xx$ and 
(b) if additionally $\gamma( N( \{ u \} | \xx\cup\{ u\}) ) > 0$ then 
$\gamma(\xx\cup\{ u\}) > 0$. 

\begin{thm} \label{t:HC-tree}
Let $(V, E, \Phi)$ be a graph with Euclidean edges such that
$(V, E)$ is a tree and $L$ the corresponding network.
Let $p$ a measurable function from the set of finite 
point configurations into $[0, \infty)$. Then $p$ is a Markov function on 
the family of configurations in general position with respect to the 
Delaunay relation if and only if
\[
p(\xx) \propto 1\{ \gamma(\xx) > 0 \} \, \prod_{x_i\in\xx} \gamma( \{ x_i \}) \,
  \prod_{i < j: x_i \sim_{\xx} x_j } \gamma(\{ x_i,  x_j \}) 
\]
for some measurable, non-negative interaction function $\gamma$.
\end{thm}

If $p$ can be normalised into a probability density, e.g. by assuming
that $\gamma( \{ x_i \})$ is bounded and $\gamma( \{ x_i, x_j \}) \leq 1$,
$p$ is a Markov density and we can define a nearest-neighbour Markov 
point process by (\ref{e:dens}). 

\

The characterisation allows some flexibility in the choice of $\gamma$.

\begin{ex}
Set $\gamma(\emptyset) = \alpha > 0$, $\gamma(\{ x_i \}) \equiv \beta > 0$,
and suppose that the interaction function for pairs depends only
on the weighted shortest path distance between its elements:
\[
\gamma(\{ x_i, x_j \}) = g( d_G( x_i,  x_j ) )
\]
for some function $g: [0,\infty)\to [0,\infty)$.
For configurations $\xx$ with $n(\xx) > 2$, set
\[
\gamma(\xx) = 
1 \{ g( d_G( x_i, x_j ) ) > 0 \mbox{ for all } x_i, x_j \in \xx \}.
\]
Clearly $\gamma(\xx) > 0$ if and only if $g( d_G( x_i, x_j ) ) > 0$ 
for all $x_i, x_j \in \xx$. 

The function $\gamma$ thus defined is an interaction function provided
that if the conditions
\begin{itemize}
\item $g( d_G( x_i, x_j ) ) > 0$ for all $x_i, x_j \in \xx$ and
\item $g( d_G( x_i, u ) ) > 0$ for all $x_i \in \xx$ such that
$x_i \sim_{\xx \cup \{ u \}} u$
\end{itemize}  
hold, then also $g( d_G( x_i, u ) )>0$ for all $x_i \in \xx$.
An alternative condition is \cite[(G) in \S 5.2]{BaddMoll89} in
combination with the interaction function
\(
\gamma(\xx) = 1 \{ g( d_G( x_i, x_j ) ) > 0 \mbox{ for all } x_i, x_j \in \xx
\mbox{ such that } x_i \sim_\xx x_j \}
\)
for configurations $\xx$ with $n(\xx) > 2$. In either case, a
sufficient condition is that the function $g$ takes strictly 
positive values. Then the Papangelou conditional intensity reads
\[
\lambda( u | \xx ) = \frac{ p(\xx\cup\{ u\} ) }{ p(\xx) } = \beta \,
\frac{ \prod_{x \sim_{\xx\cup \{ u\} } u} g( d_G( x, u) ) }{
       \prod_{x \sim_{\xx} y, x\not\sim_{\xx\cup \{ u\}} y} g( d_G(x, y) ) },
\]
in complete analogy with the conditional intensity of a renewal process
on the real line.
\end{ex}

\subsection{The local Delaunay relation}

As shown in Example~\ref{e:triangle}, the Delaunay relation does
not satisfy the consistency relations (C1)--(C2) if the graph
$(V,E)$ is not a tree. However, we may employ the neighbourhood
relation implicit in the graph to define a local Delaunay relation.
Such a procedure is similar to one employed in image analysis for
edge detection and texture analysis \cite{Gemaetal90}.

\begin{defn} \label{d:edge}
Let $(V, E, \Phi)$ be a graph with Euclidean edges and $L$ its
associated network. Define a symmetric and reflexive relation $\sim_E$ 
on $L$ as follows:
\[
(i,x) \sim_E (j,y)  \Leftrightarrow \left\{ \begin{array}{ll}
\phi_i^{-1}(\partial J_i) \cap \phi_j^{-1}(\partial J_j) \neq \emptyset,
&  i, j \neq 0 \\
\phi_i^{-1}(\partial J_i) \cap \{ y \} \neq \emptyset,
&  i \neq 0, j = 0 \\
\{ x, y \} \in E, 
&  i = j = 0.
\end{array} \right.
\]
Write, for $i,j \neq 0$, $\sim_{\xx}^{i,j}$ for the Delaunay 
relation restricted to
\[
L \cap \left( ( \{ i, 0 \} \times \phi_i^{-1}(\bar J_i) ) 
\cup 
( \{ j, 0 \} \times \phi_j^{-1}(\bar J_j) )  \right),
\]
the restriction of $L$ to at most two edges and their incident
vertices, and define
\[
(i,x) \sim^E_\zz (j,y) \Leftrightarrow \left\{ \begin{array}{l}
(i,x) 
\sim^{i,j}_{\zz \cap ( \{ i, 0 \} \times \phi_i^{-1}(\bar J_i)) 
\cup (\{ j,0 \} \times \phi_j^{-1}(\bar J_j)) } 
(j,y);
(i,x) \sim_E (j,y), \\
\quad \quad \quad  i, j \neq 0 \\
(i,x) 
\sim^{i,i}_{\zz \cap ( \{ i, 0 \} \times \phi_i^{-1}(\bar J_i)) } 
(0,y);
(i,x) \sim_E (0,y), \\
\quad \quad \quad i \neq 0, j = 0 \\
(0,x) \sim^{k,k}_{ \zz \cap (\{ k, 0 \} \times \phi_k^{-1}(\bar J_k)) } (0,y);
 \{ x, y \} = e_k \in E, \\
\quad \quad \quad  i = j = 0.
\end{array} \right.
\]
\end{defn}

In words, a vertex is a $\sim_E$-neighbour of the edges it is incident with 
and edges are neighbours if they share a common vertex. Thus, the 
$\sim_E$-relation does not depend on the configuration. It does, however, 
crucially depend on the geometry of the graph -- splitting an edge will 
result in a different relation. After combination with the Delaunay relation,
the resulting relation $\sim^E_\zz$ is configuration dependent and depends
on the geometry of the graph.

The main result of this section is the following. Its proof will be
given in Section~\ref{S:theorems}.

\begin{thm} \label{t:local}
Let $(V, E, \Phi)$ be a graph with Euclidean edges and $L$ its associated 
network. Then the local Delaunay relation $\sim_{\xx}^E$ satisfies (C1)--(C2).
\end{thm}

\section{Proofs of Lemmata}
\label{S:proofs}

\begin{proof}{Lemma~\ref{L:hereditary}}\\
We claim that $y_1 \sim_{\xx} y_2$ implies that $y_1 \sim_{\zz} y_2$
for any $\zz\subseteq \xx$ and $y_1, y_2 \in \zz$. To see this, note that 
\[
\chi(\{ y_1, y_2 \} | \xx) = 1 \Leftrightarrow \exists \xi \in L:
d_G(\xi, x) \geq d_G(\xi, y_i) = d_G(\xi, y_2) \mbox{ for all } x \in \xx.
\]
A fortiori, for this $\xi$,
\(
d_G(\xi, z) \geq d_G(\xi, y_i) = d_G(\xi, y_2) 
\)
for all $z$ in the smaller set $\zz \subseteq \xx$, 
so that $\chi(\{ y_1, y_2 \}) | \yy) = 1$, that is, $y_1 \sim_{\zz} y_2$.

By convention, the clique indicator function takes the value one for 
singletons and the empty set regardless of the configuration. Hence
we may focus on configurations $\yy$ of cardinality at least two.

For the first statement, suppose that $\chi(\yy|\xx) = 1$. Pick any 
$y_1, y_2 \in \yy$. Then $y_1 \sim_\xx y_2$ and, by the above claim, 
$y_1 \sim_\zz y_2$. Since $y_1$ and $y_2$ are chosen arbitrarily, 
$\chi(\yy|\zz) = 1$.  

For the second statement, if $\chi(\yy|\zz) = 0$, there exist
$y_1, y_2 \in \zz$ such that $y_1 \not\sim_\zz y_2$. By the claim,
also $y_1 \not \sim_\xx y_2$, hence $\chi(\yy|\xx) = 0$.
\end{proof}

\begin{proof}{Lemma~\ref{L:unique}}\\
Suppose that $(V,E,\Phi)$ is a tree. For vertices, since $(V,E)$ is
a tree, there is a unique path between each pair of vertices. Hence
we may restrict ourselves to a pair of points $(i,x)$ and $(j,y)$ of
which at least one belongs to the set $L \setminus (\{ 0 \} \times V)$.

If $i=j$, one path runs along the edge. Since a walk does not reverse
its tracks by definition, any other walk from $(i,x)$ to $(i,y)$ must
visit at least one of the two end vertices of $e_i$ before returning,
hence use $e_i$ twice. 

If $i=0$ and $j\neq 0$, note there is a unique path from $x$ to 
$v_j^1$ in the graph $(V,E)$ and hence a corresponding one in the 
labelled space $L$. If this path includes edge $e_j$, by deleting the
$j$-labelled curve segment from $y$ to $v_j^1$, one obtains a path 
from $(i,x)$ to $(j,y)$,  otherwise such a path is found by extending
the path to $v_j^1$ with this segment. Moreover, any other path would 
have to pass at least one of the vertices in $e_j$, that is, coincide 
with the construction up to there and hence entirely. 

Finally, if $i,j\neq 0$, as seen in the previous case, there is 
a unique path from $(0,v_i^1)$ to $(j,y)$. If this path includes edge $e_i$,
this yields the path from $(i,x)$ to $(j,y)$. Otherwise, precede by the
segment along $e_i$ from $(i,x)$ to $(0,v_i^1)$. Any other path would have
to pass at least one of the vertices in $e_i$, that is, coincide with
the construction from there and hence entirely.

Reversely, let $(V, E, \Phi)$ be such that there is a unique path between
any pair of distinct points $(i,x)$ and $(j,y)$ in $L$ (which exists
since a graph with Euclidean edges is connected by definition). 
In particular, there is a unique path between any pair of $0$-labelled
vertices, and therefore $(V, E)$ is a tree.
\end{proof}

\begin{proof}{Lemma~\ref{L:wheel}} \\
Since $L$ is a tree, by Lemma~\ref{L:unique}, 
there are unique paths from $y_1$ to $y_2$ and 
from $y_2$ to $y_3$, say with consecutive labelled vertices
\[
y_1, (0, v_1), \dots, (0,v_p), y_2, (0,v_{p+1}), \dots, (0,v_{p+q}), y_3,
\]
$p, q \geq 0$ and no vertex replication in $v_1, \dots, v_p$ nor in 
$v_{p+1}, \dots, v_{p+q}$.

Suppose 
that $v_k = v_l$ for some $1\leq k \leq p$ and $p+1 \leq l \leq p+q$. 
Then the paths 
\(
((0,v_k), \dots, (0,v_p), y_2)
\)
and
\(
((0,v_l), \dots, (0,v_{p+1}), y_2)
\)
both connect $(0, v_k) = (0, v_l)$ and $y_2$ in $L$, and must therefore 
coincide. Extending if possible, we may assume that $v_{k-1} \neq v_{l+1}$.
If some earlier vertex $v_i$, $i<k-1$ would be identical to $v_j$ for
$l+1 < j \leq p+q$, a cycle would be created from $v_i$ via $v_k = v_l$ 
to $v_j = v_i$. Hence the sequences $(v_1, \dots, v_{k-1})$ and 
$(v_{l+1}, \dots, v_{p+q})$ do not intersect and the paths 
$(y_1, (0, v_1), \dots, (0, v_{k-1}))$, 
$(y_2, (0, v_{p}), \dots, (0, v_{k+1}))$ and
$(y_3, (0,v_{p+q}), \dots, (0, v_{l+1}))$ 
are connected at $(0, v_k) = (0, v_l )$ and therefore form three stokes
connected at a single hub provided the lengths are positive. It therefore
remains to consider the cases $k=1$, $k=p$ or $l=p+q$.

If $l=p+q$, $y_3$ may lie on the path from $y_1$ to $y_2$; if it does
not, the curve segment from $(0,v_l)$ to $y_3$ forms a stoke of 
positive length. Similarly, if $k=p$, $y_2$ may lie on the path from
$y_1$ to $y_3$; if it does not it, the curve segment from $(0, v_k)$ to
$y_2$ forms a stoke of positive length. Also, if $k=1$, $y_1$ may
lie along the path between $y_2$ and  $y_3$; if it does not, the
curve segment from $y_1$ to $(0, v_k)$ forms a stoke of positive length.

Finally, if there is no $1 \leq k \leq p$ such that $v_k = v_l$ for some 
$p+1 \leq l \leq p+q$, the unique path from $y_1$ to $y_3$ runs via $y_2$.
\end{proof}

\begin{proof}{Lemma~\ref{L:clique}} \\
Suppose that $\chi(\yy|\xx) = 1$ for some $\yy = \{ y_1, y_2, y_3 \}
\subseteq \xx \subseteq L$. By Lemma~\ref{L:hereditary}, $\chi(\yy|\yy)
= 1$ hence one needs only consider $\xx = \yy$.

By Lemma~\ref{L:wheel}, the elements of $\yy$ either form a wheel or
a path. First consider the case that $\yy$ is a wheel with three
stokes emanating from a hub $(0,v)$. Ordering the stokes according to
their length, without loss of generality suppose that 
$a\leq b \leq c$ where $a = d_G( y_1, (0, v) )$, $b = d_G( y_2, (0, v))$ 
and $c = d_G( y_3, (0, v))$. Since by assumption $\yy$ is in general 
position, $a<b<c$. Therefore, as $L$ is a tree, 
$C(y_2|\yy) \cap C(y_3| \yy) = \emptyset$ and $y_2 \not \sim_\yy y_3$. 

It remains to consider the case that all three elements of $\yy$ lie
along a path, without loss of generality from $y_1$ via $y_2$ to $y_3$. 
Then $C(y_1|\yy) \cap C(y_3|\yy) = \emptyset$, again using the 
tree property to ensure that there are no paths to connect $y_1$ and
$y_3$ other than via $y_2$. Therefore $y_1 \not \sim_\yy y_3$.
In conclusion, there cannot be a clique of cardinality three or larger.

\bigskip

To prove the second assertion, let $\yy = \{ y_1, y_2\}$ and write
$\xi$ for the midpoint with respect to $d_G$ along the path between
$y_1$ and $y_2$. Since $L$ is assumed to be a tree, by Lemma~\ref{L:unique},
the path is unique. Since all $\phi_i \in \Phi$ are bijections,
also the midpoint is unique. Let $\xx \supseteq \yy$ be some 
configuration in general position.

If $\xi \in C(y_1|\xx) \cap C(y_2|\xx)$, clearly $\chi(\yy|\xx) = 1$.
Reversely, suppose that $\yy$ is a clique in $\xx$. Then there exists
some $\eta\in L$ such that $\eta \in C(y_1 | \xx ) \cap C(y_2 | \xx)$.
By Lemma~\ref{L:wheel}, the triple $\{ \eta, y_1, y_2 \}$ either forms
a path or a wheel. In the first case, the property of equidistance to 
$y_1$ and $y_2$ implies that $\eta = \xi$ and the proof is complete.

Next consider the case that $\eta$ and $\yy$ form three stokes of
strictly positive length connected at some hub $(0,v)$. Since 
\(
d_G(\eta, y_j) = d_G(\eta, (0,v)) + d_G((0,v), y_j),
\)
$j=1,2$, by the uniqueness of the paths from $\eta$ to $y_j$ 
(cf.~Lemma~\ref{L:unique}), $d_G( (0,v), y_1) = d_G( (0,v), y_2)$.
Hence $\xi = (0,v)$, now using the uniqueness of the path between $y_1$ 
and $y_2$. It remains to show that no other point of $\xx$ lies
closer to $\xi$ than the $y_j$, $j=1, 2$. Since $L$ is a tree, any 
such $x\in\xx$ is connected either to the hub $\xi$ or to exactly 
one of the tree stokes. Any $x\in\xx$ connected to the end points 
$y_1, y_2$ of the stokes lies further from $\xi$ than $y_1$ and $y_2$. 

Any $x\in\xx$ connected to the stoke of $\eta$ satisfies,
since $\eta$ is assumed to lie in $C(y_1 | \xx ) \cap C(y_2 | \xx)$,
\[
d_G(\eta, x) \geq d_G(\eta, y_j) = d_G(\eta, \xi) + d_G(\xi, y_j), 
   \quad j = 1, 2,
\]
again using the uniqueness of paths for the last equality.
If the connection is at the endpoint $\eta$, then
\(
d_G(\xi, x) = d_G(\xi, \eta) + d_G(\eta, x)
\)
hence, by the above equation,
\(
d_G(\xi, x) \geq 2 d_G(\xi, \eta) + d_G(\xi, y_j) \geq d_G(\xi, y_j).
\)
If the connection is at a vertex $(0,w)$, $w\in V$,
\begin{eqnarray*}
d_G(\xi, x) & = & d_G(\xi, (0,w)) + d_G(x, (0,w)) \\
& = & d_G(\xi, (0,w)) + d_G(\eta, x) - d_G(\eta, (0,w) )
\\
& \geq &
d_G(\xi, (0,w) ) + d_G(\eta, y_j) - d_G(\eta, (0,w) )\\
& = &
d_G(\xi, (0,w) ) + d_G( y_j, (0,w)) \geq d_G(\xi, y_j).
\end{eqnarray*}

Finally, any $x\in\xx$ connected to the stoke of a $y_j$ at
$\xi$ or some other vertex $(0,w)$ satisfies, since $\eta$ is
assumed to lie in $C(y_1 | \xx ) \cap C(y_2 | \xx)$,
\[
d_G(\eta, x) \geq d_G(\eta, y_j) = d_G(\eta, \xi) + d_G(\xi, y_j)
\]
so that
\(
d_G(\xi, x) = d_G(\eta, x) - d_G(\eta, \xi)  \geq
d_G(\xi, y_j)
\)
and the proof is complete.
\end{proof}

\section{Proof of main theorems}
\label{S:theorems}

\begin{proof}{Theorem~\ref{t:Delaunay}} \\
First, note that for $\yy$ of size zero or one, (C1) and (C2) are 
automatically satisfied. The assumption that $L$ is a tree 
and the restriction to configurations in general position imply, by 
Lemma~\ref{L:clique}, that $\chi(\yy| \zz) = 0$ for all $\zz$ 
when the cardinality of $\yy$ is three or more. Hence (C1) and 
(C2) hold for $\yy$ with $n(\yy) \geq 3$ as well and it suffices 
to consider pairs $\yy = \{ y_1, y_2 \}$.

\paragraph{Condition (C1)}
First, consider (C1). 
Take $y_1, y_2 \in \zz \subseteq L$, $u \in L$ with $u \not \in \zz$,
and suppose that $\chi( \{ y_1, y_2 \} | \zz )$ differs from 
$\chi( \{ y_1, y_2 \} | \zz \cup \{ u \} )$.  By Lemma~\ref{L:hereditary}, 
if $\chi( \{ y_1, y_2 \} | \zz \cup \{ u \}) = 1$, also
$\chi( \{ y_1, y_2 \} | \zz ) = 1$, so it suffices to consider the 
case that $\chi( \{ y_1, y_2 \} | \zz \cup \{ u \}) = 0$ but
$\chi( \{ y_1, y_2 \} | \zz ) = 1$. 

Let $\tilde y_j$ be the point lying halfway between $y_j$ and $u$
along the unique path between them (cf.~Lemma~\ref{L:unique}) and 
let $\xi$ be the halfway point between $y_1$ and $y_2$. These points 
exist since the parametrisations are bijections. 
By Lemma~\ref{L:clique}, since $\chi(\yy|\zz) = 1$, $\xi \in
C(y_1 | \zz ) \cap C(y_2 | \zz )$, so that for all $z\in \zz$,
$d_G(\xi, z) \geq d_G(\xi, y_j)$, $j=1,2$. Also by construction 
\(
d_G(\tilde y_j, y_j) = d_G(\tilde y_j, u).
\)
We shall show that, for $j=1, 2$,
\(
\tilde y_j \in C(y_j | \zz \cup \{ u \} ) \cap C(u | \zz \cup \{ u \})
\)
and therefore $y_j \sim_{\zz\cup \{u\}} u$. To do so, we
must demonstrate that 
\begin{equation}
d_G( \tilde y_j, z ) \geq d_G(\tilde y_j, u) = d_G(\tilde y_j, y_j)
\label{e:tmp-z}
\end{equation}
for $z\in \zz \setminus \{ y_j \}$.

By Lemma~\ref{L:wheel}, since $L$ is a tree, the paths between the 
three points $u, y_1, y_2$ form either a wheel with three stokes 
joined at some hub $\eta = (0,v)$ with $v\in V$, or the points lie 
on a path.

First, consider the case that $u$, $y_1$ and $y_2$ lie on a path.  
If $u$ were an extremity of this path, say $y_1$ is on the path 
between $u$ and $y_2$, then
\[
d_G(\xi, u ) = d_G(\xi, y_1) + d_G(y_1, u) > d_G(\xi, y_1).
\]
By the assumption $\chi(\yy|\zz) = 1$ and Lemma~\ref{L:clique}, 
for any  $z\in \zz$, also $d_G(\xi, z) \geq d_G(\xi, y_1) = d_G(\xi, y_2)$ 
so that $\xi \in C(y_1 | \zz \cup \{ u \} ) \cap C(y_2 | \zz \cup \{ u \})$,
thus violating the assumption that $\chi(\yy|\zz \cup \{ u \}) = 0$.
We conclude that $u$ has to lie on the path from $y_1$ to $y_2$. 
Suppose that (\ref{e:tmp-z}) does not hold, i.e. that for some $j$ and
some $z \in \zz \setminus \{ y_j \}$,
the distance $d_G(\tilde y_j, z) < d_G(\tilde y_j, y_j)$.
Then
\[
d_G(\xi, z) \leq d_G(\xi, \tilde y_j) + d_G(\tilde y_j, z) <
d_G(\xi, \tilde y_j) + d_G(\tilde y_j, y_j) = d_G(\xi, y_j)
\]
using uniqueness of paths (Lemma~\ref{L:unique}).
However, if $d_G(\xi, z) < d_G(\xi, y_j)$, then $\xi \not \in C(y_j|\zz)$ 
in contradiction with Lemma~\ref{L:clique}.

Next suppose that $\{ u, y_1, y_2 \}$ form a wheel and, without loss of 
generality, the path from 
$y_1$ to $y_2$ passes first $\xi$ and then $\eta$. Since $\chi(\yy| \zz \cup \{ u \})
= 0$, the intersection of the $C(y_j|\zz\cup \{ u \})$, $j=1,2$,
is empty and in particular does not contain $\xi$. Hence 
\[
d_G(\xi, u) < d_G(\xi, y_1) = d_G(\xi, y_2). 
\]
Therefore, 
\[
d_G(\eta, u) < d_G(\eta, y_2) \leq d_G(\eta, y_1)
\]
with equality only if $\eta = \xi$. Consequently, $\tilde y_j$ lies
on the stoke of $y_j$ for each $j=1, 2$. 
Now, if $d_G(\tilde y_2, z) < d_G(\tilde y_2, y_2)$ for some $z\in \zz 
\setminus \{ y_2 \}$, since
$\tilde y_2$ lies on the path from $\xi$ to $y_2$ via $\eta$, 
then
\[
d_G(\xi, z) \leq d_G(\xi, \tilde y_2) + d_G(\tilde y_2, z) <
d_G(\xi, \tilde y_2) +  d_G(\tilde y_2, y_2) = d_G(\xi, y_2)
\]
in contradiction with the assumption that $\xi \in C(y_2 | \zz)$.
Similarly for $y_1$, recalling that $d_G(\xi, u) < d_G(\xi, y_1)$, 
the point $\tilde y_1$ lies on the path between $\xi$ and $y_1$.
Hence if $d_G(\tilde y_1, z) < d_G(\tilde y_1, y_1)$, then
\[
d_G(\xi, z) \leq d_G(\xi, \tilde y_1) + d_G(\tilde y_1, z) <
d_G(\xi, \tilde y_1) +  d_G(\tilde y_1, y_1) = d_G(\xi, y_1)
\]
in contradiction with the assumption that $\xi \in C(y_1 | \zz)$.
Therefore (\ref{e:tmp-z}) cannot be violated and (C1) holds.

\paragraph{Condition (C2)}
Next, consider (C2). Take $y_1, y_2 \in \zz \subseteq L$ and $u,v \in L$ 
such that $u, v \not \in \zz$ with $u \not \sim_\xx v$ where $\xx =
\zz \cup \{ u, v \}$. Write $\yy = \{ y_1, y_2 \}$. 
Lemma~\ref{L:hereditary} implies that if
$\chi(\yy | \zz) = 0$, the same is true upon adding points to $\zz$ and
(C2) holds in this case. Therefore, it suffices to consider the case that
$\chi( \{ y_1, y_2 \} | \zz ) = 1$. The same lemma implies that if
$\chi( \{ y_1, y_2 \} | \zz \cup \{ u, v \} ) = 1$, this remains true
when deleting points from $\zz$ and (C2) holds. The only case left to
consider is that when $\chi( \{ y_1, y_2 \} | \zz ) = 1$ and
$\chi( \{ y_1, y_2 \} | \zz \cup \{ u, v \} ) = 0$. We must show that
exactly one of $\chi( \{ y_1, y_2 \} | \zz \cup \{ u \})$ and
$\chi( \{ y_1, y_2 \} | \zz \cup \{ v \})$ takes the value $1$ and will
do so by contradiction.

Let $\eta\in L$ be the point that lies halfway between $u$ and $v$ along the
unique path between them (cf.\ Lemma~\ref{L:unique}) and write $\xi$ 
for the halfway point between $y_1$ and $y_2$. These points exist
since the parametrisations are bijections. 

Suppose that $\chi(\yy|\zz\cup \{ u \} ) = 0 = \chi(\yy|\zz\cup \{ v \} )$. 
By the proof of condition (C1) above, the triple $\{ u, y_1, y_2 \}$ form 
a direct path with $u$ on the path from $y_1$ to $y_2$, or a wheel with
three stokes connecting in a hub. In either case 
\[
d_G(\xi, u) < d_G(\xi, y_1) = d_G(\xi, y_2).
\]
The same is true for the triple $\{ v, y_1, y_2 \}$. Therefore the ensemble 
can be seen as a path from $y_1$ to $y_2$ that passes branches leading off
to $u$ and to $v$ if we include the degenerate cases of a branch consisting 
of the single point $u$, respectively $v$.
Without loss of generality, suppose the order is $y_1$, then the stokes of
$u$ and $v$, and finally $y_2$. We claim that such an arrangement would 
imply that $u \sim_{\xx} v$ in contradiction with the assumption.

To prove the claim, we show that $\eta \in C(u | \zz \cup \{ u, v \})
\cap C(v| \zz \cup \{ u, v \})$, that is $d_G(\eta, u) = d_G(\eta, v)
\leq d_G(\eta, z)$ for all $z\in \zz$. Suppose not. Then for some 
$z\in \zz$, possibly $y_1$ or $y_2$, $d_G(\eta, z) < d_G(\eta, u) =
d_G(\eta, v)$ and therefore
\[
d_G(\xi, z) \leq d_G(\xi, \eta) + d_G(\eta, z) <
d_G(\xi, \eta) + d_G(\eta, u)  = d_G(\xi, \eta) + d_G(\eta, v).
\]
The right hand side is equal to either $d_G(\xi, u)$ or $d_G(\xi, v)$
and recalling that both are strictly smaller than $d_G(\xi, y_j)$
we obtain
\(
d_G(\xi, z) < d_G(\xi, y_1) = d_G(\xi, y_2)
\)
in contradiction with the assumption that $\chi(\yy|\zz) = 1$,
that is, $\xi \in C(y_j|\zz)$ for $j=1,2$ (Lemma~\ref{L:clique}). Hence 
$C(u | \zz \cup \{ u, v \}) \cap C(v| \zz \cup \{ u, v \})$ contains
$\eta$, implying $u \sim_{\xx} v$ in contradiction with the assumption.

Finally, suppose that $\chi(\yy|\zz\cup \{ u\}) = 1 = \chi(\yy|\zz\cup\{v\})$.
By Lemma~\ref{L:clique}, 
$\xi \in C(y_1 | \zz \cup \{ u \} ) \cap C(y_2 | \zz \cup \{ u \} )$ 
and therefore
$d_G(\xi, u) \geq d_G(\xi, y_1) = d_G(\xi, y_2)$. Similarly, 
\(
d_G(\xi, v) \geq d_G(\xi, y_1) = d_G(\xi, y_2).
\) 
By assumption, $\chi(\yy | \zz \cup \{ u, v \}) = 0$, which 
means that $C(y_1 | \zz \cup \{ u, v \} ) \cap C(y_2 | \zz \cup \{ u, v\}) 
= \emptyset$ and in particular does not contain $\xi$. Therefore,
recalling the assumption $\chi(\yy | \zz) = 1$, we have that
$\min( d_G(\xi, u), d_G(\xi, v) ) < d_G(\xi, y_1) = d_G(\xi, y_2)$,
a contradiction.

In conclusion, exactly one of $\chi(\yy|\zz\cup\{u\})$ and 
$\chi(\yy|\zz\cup\{v\})$ takes the value one.
\end{proof}

\begin{proof}{Theorem~\ref{t:HC-tree}}\\
Suppose that $p$ is a Markov function. By Theorem~\ref{t:Delaunay}, 
the Delaunay relation satisfies (C1) and (C2) and therefore 
\cite[Thm~4.13]{BaddMoll89} implies that $p$ can be factorised as
\[
p(\xx) = \prod_{\yy \subseteq \xx} \gamma(\yy)^{\chi(\yy|\xx)} 
\]
for some interaction function $\gamma$ under the convention that
$0^0 = 0$. Since by Lemma~\ref{L:clique}
the cliques have cardinality at most two, the factorisation reduces to
\[
p(\xx) \propto \prod_{i} \gamma( \{ x_i \}) 
  \prod_{i < j} \gamma(\{ x_i,  x_j\} )^{\chi(\{x_i, x_j\}|\xx)} 
  \left( \prod_{\yy\subseteq \xx: n(\yy) > 2} \gamma(\yy) \right)^0.
\]
As $\gamma$ is an interaction function, it is hereditary, so 
$\gamma(\xx) > 0$ if and only if $\gamma(\yy) > 0$ for all $\yy\subseteq \xx$,
hence 
\[
\left( \prod_{\yy\subseteq \xx: n(\yy) > 2} \gamma(\yy) \right)^0 =
1\{ \gamma(\xx) > 0 \}.
\]
When $\gamma(\xx) > 0$, also $\gamma(\{ x_i \})$ and 
$\gamma(\{ x_i, x_j\})$ are strictly positive, so
\[
p(\xx) = \gamma(\emptyset) \prod_i \gamma( \{ x_i \} ) 
      \prod_{i<j, x_i \sim_\xx x_j} \gamma( \{ x_i, x_j \} ).
\]
Clearly $p(\xx) = 0$ when $\gamma(\xx) = 0$.

Reversely, any function of the specified form is a Markov 
function. To see this, suppose $p(\xx) > 0$. Then $\gamma(\xx) > 0$ and
for all subsets $\yy$ of $\xx$ also $\gamma(\yy) > 0$. Hence $p(\yy) > 0$.
Furthermore, the ratio
\(
{p(\xx\cup \xi)} / {p(\xx)}
\)
can be written as
\[
\gamma(\xi) \, 1\{ \gamma(\xx\cup\{ \xi \}) > 0 \} \,
\prod_{i} \gamma( \{ x_i, \xi \})^{\chi( \{ x_i, \xi \} | \xx \cup \{ \xi \} ) }
\times
\prod_{i<j} \gamma( \{ x_i, x_j\} )^{\chi( \{ x_i, x_j \} | \xx \cup \{ \xi \} ) 
- \chi( \{ x_i, x_j \} | \xx ) }.
\]
By Theorem~\ref{t:Delaunay}, the Delaunay relation satisfies (C1), 
hence the product
\[
\prod_{i<j} \gamma( \{ x_i, x_j\} )^{\chi( \{ x_i, x_j \} | \xx \cup \{ \xi \} ) 
- \chi( \{ x_i, x_j \} | \xx ) }
\]
depends only on $x_i$ and $x_j$ in $N(\{ \xi \} | \xx \cup \{ \xi \} )$
and the relations $\sim_\xx$ and $\sim_{\xx\cup\{ \xi\} }$ on this neighbourhood. 
Furthermore, recalling that $\gamma$ is an interaction function and
$\gamma(\xx) > 0$, 
\[
 1\{ \gamma(\xx\cup\{ \xi \}) > 0 \}
\prod_{i} \gamma( \{ x_i, \xi \})^{\chi( \{ x_i, \xi \} | \xx \cup \{ \xi \} ) }
= 1\{ \gamma(N( \{ \xi \} | \xx\cup\{ \xi \})) > 0 \}
\prod_{\xi \sim_{\xx\cup\{ \xi\}} x_i} \gamma( \{ x_i, \xi \})
\]
depends only on $N(\{ \xi \} | \xx\cup\{ \xi \}))$ and the relation 
$\sim_{\xx\cup \{ \xi\} }$ on this neighbourhood.
\end{proof}

\begin{proof}{Theorem~\ref{t:local}}
Let $\yy \subseteq \zz \subseteq L$ and $u,v\in L$ be such that 
$u,v\not \in \zz$ and $u \not \sim_{\zz\cup\{ u, v\} } v$.
 
We first observe that cliques in $\sim_E$ consist of points lying on 
closed edges emanating from a single vertex. Thus, the size of a 
$\sim_E$-clique may be larger than two. Moreover, the  
Delaunay relation restricted to a pair of such edges coincides with the 
sequential neighbourhood relation with respect to an ordering defined by the 
parametrisations on the edges. In other words, consecutive points on a 
single edge are each other's nearest neighbours; also the point on one of
the edges that is closest to the vertex that joins the two edges, if it 
exists, is a nearest neighbour of the point closest to that vertex on the 
other edge, if it exists, and no other points are nearest neighbours.
In particular, each point has at most two nearest neighbours. Cliques
in the combined relation $\sim^E_\zz$ therefore are either empty, consist 
of a single point, of two consecutive points on a single edge, or
of points on different edges that are closest on their edge to the 
central vertex from which all edges emanate. The clique size is 
therefore at most the degree of the central vertex.

\paragraph{Condition (C1)}
If $\chi(\yy|\zz) = 0$ there exists a pair $y_1 \neq y_2 \in \yy$
for which $y_1 \not \sim^E_\zz y_2$. Then either $y_1$ and $y_2$ lie
on edges that are not adjacent, in which case 
$y_1 \not \sim^E_{\zz \cup \{ u \}} y_2$ and $\chi(\yy|\zz\cup \{ u \}) = 0$,
or $y_1$ and $y_2$ lie on adjacent edges but $y_1$ and $y_2$ are not
sequential neighbours. The addition of $u$ cannot make them sequential
neighbours, so $\chi(\yy|\zz\cup\{ u \}) = 0$. 

Suppose therefore that $\chi(\yy|\zz) = 1$ but there exists a pair 
$y_1 \neq y_2 \in \yy$ for which $y_1 \not \sim^E_{\zz\cup \{ u \}} y_2$.
Then $y_1$ and $y_2$ must lie on $\sim_E$-related edges (either a single
one, or two adjacent ones), be consecutive in $\zz$ but not in $\zz
\cup \{ u \}$. This can only happen if $u$ lies in between $y_1$ and
$y_2$, making $y_1$ and $y_2$ both $\sim^E_{\zz\cup \{ u \} }$-neighbours
of $u$. 

\paragraph{Condition (C2)}
As shown when proving (C1), if $\chi(\yy|\zz) = 0$, also $\yy$ cannot
be a clique in configurations with more points. If 
$\chi(\yy| \zz \cup \{ u, v\} ) = 1$, all points in $\yy$ lie on
adjacent edges and are sequential neighbours in the set $\zz \cup \{
u, v\})$ restricted to their edges. The same remains true when 
$u$ and $v$ are removed, so that $\chi(\yy | \zz) = 1$. 
Hence it remains to consider the case that $\chi(\yy|\zz) = 1$ but 
$\chi(\yy|\zz\cup\{ u, v \}) = 0$. In this case, the points of $\yy$ 
must lie on a single edge or on a number of edges that emanate from 
a single vertex $(0,w)\in L$, $w\in V$. 

Now, if $u$ and $v$ do not lie on any of these $\yy$-edges, 
$\chi(\yy|\zz\cup \{ u, v\}) = 1$, in contradiction with the assumptions. 

If exactly one of $u$ and $v$ does not lie on any of the $\yy$-edges, 
without loss of generality $u$, then 
$\chi(\yy|\zz\cup\{ u \}) = \chi(\yy|\zz) = 1$. Moreover, since by 
assumption $\chi(\yy|\zz\cup \{ u, v \}) = 0$, there is a pair 
$y_1, y_2\in \yy$ that are not consecutive in the configuration 
$\zz \cup \{ u, v\}$ restricted to the edge or edges on which 
$y_1$ and $y_2$ lie. Since they are adjacent in the configuration $\zz$, 
$v$ must lie in between $y_1$ and $y_2$, which implies that 
$\chi(\yy| \zz \cup \{ v \}) = 0$ in accordance with (C2).

Finally, suppose both $u$ and $v$ lie on the $\yy$-edges that emanate
from $(0, w)$. If $\yy = \{ y_1, y_2 \}$ consists of two consecutive 
points in $\zz$, 
since $y_1$ and $y_2$ are no longer consecutive in the configuration 
$\zz \cup \{ u, v \}$, they must be separated by either $u$ or $v$
but not by both, since by assumption $u$ and $v$ are not sequential
neighbours in $\zz \cup \{ u, v\}$. Hence (C2) holds. If $\yy$ 
consists of more than two points $\yy = \{ y_1, \dots, y_k \}$, the 
assumption $\chi(\yy|\zz) = 1$ implies that the $y_j$, $j=1, \dots, k$, 
must lie on different edges $e_1, \dots, e_k$ emanating from $w$ and 
no points of $\zz$ lie between $y_i$ and $w$. Since $u$ and $v$ are 
not sequential neighbours, they cannot both lie between some $y_i$ 
and $(0, w)$; one of them, however, must, since the clique indicator 
function of $\yy$ in $\zz\cup \{ u , v\}$ takes the value zero. 
Consequently, also in this case (C2) is seen to hold.
\end{proof}


\begin{thebibliography}{99}
\bibitem{Moll16}
Anderes, E., M\o ller, J. and Rasmussen, J.G. (2016).
Second-order pseudo-stationary random fields and point processes
on graphs and their edges.
Talk presented by J.\ M\o ller at {\em AU Workshop on Stochastic 
Geometry, Stereology and their Applications\/}, Sandbjerg Manor.
Extended and generalised in ArXiv manuscript 1710.01295.
\bibitem{BaddMoll89}
Baddeley, A.J. and M{\o }ller, J. (1989).
Nearest-neighbour {M}arkov point processes and random sets.
{\em International Statistical Review\/} {\bf 57}, 89--121.
\bibitem{Baddetal96}
Baddeley, A.J., Lieshout, M.N.M.~van and M{\o }ller, J. (1996).
Markov properties of cluster processes.
{\em Advances in Applied Probability\/} {\bf 28}, 346--355.
\bibitem{Baddetal17}
Baddeley, A., Nair, G., Rakshit, S. and McSwiggan, G. (2017).
``Stationary'' point processes are uncommon on linear networks.
{\em Stat\/} {\bf 6}, 68--78.
\bibitem{BondMurt08}
Bondy, A. and Murty, M.R. (2008).
{\em Graph theory\/}. Springer.
\bibitem{DaleVere88}
Daley, D.J. and Vere--Jones, D. (2003, 2008).
{\em An introduction to the theory of point processes\/},
second edition in two volumes. Springer.
\bibitem{Gemaetal90}
Geman, D., Geman, S., Graffigne, C. and Dong, P. (1990).
Boundary detection by constrained optimization.
{\em IEEE Transactions in Pattern Analysis Machine Intelligence\/} 
{\bf 12}, 609--628.
\bibitem{HaggLiesMoll99}
{H\"aggstr\"om}, O., Lieshout, M.N.M.~van and {M\o ller}, J. (1999).
Characterisation and simulation results including exact simulation
 for some spatial point processes.
{\em Bernoulli\/} {\bf 5}, 641--659.
\bibitem{Lies00}
Lieshout, M.N.M.~van (2000).
{\em Markov point processes and their applications\/}.
Imperial College Press.
\bibitem{OkabSugi12}
Okabe, A. and Sugihara, K. (2012).
{\em Spatial analysis along networks. Statistical and computational
methods\/}. Wiley.
\bibitem{Raksetal17}
Rakshit, S., Nair, G. and Baddeley, A. (2017).
Second-order analysis of point patterns on a network using any distance 
metric.
{\em Spatial Statistics\/}, to appear.
\bibitem{RiplKell77}
Ripley, B.D. and Kelly, F.P. (1977). 
Markov point processes. 
{\em Journal of the London Mathematical Society\/} {\bf 15}, 188--192.
\end{thebibliography}
\end{document}